\documentclass[12pt,reqno]{amsart}
\usepackage{enumerate, latexsym, amsmath, amsfonts, amssymb, amsthm, color}
\def\pmod #1{\ ({\rm{mod}}\ #1)}
\def\Z{\Bbb Z}
\def\N{\Bbb N}

\def\l{\left}
\def\r{\right}
\def\bg{\bigg}
\def\({\bg(}
\def\){\bg)}
\def\t{\text}
\def\f{\frac}

\def\ls{\leqslant}
\def\gs{\geqslant}

\def\bi{\binom}

\def\eq{\equiv}

\def\FF#1#2#3{{}_2F_1\bigg(\bmatrix{#1}\\{#2}\end{bmatrix}\bigg|#3\bigg)}
\def\Proof{\noindent{\it Proof}}

\theoremstyle{plain}
\newtheorem{theorem}{Theorem}

\newtheorem{lemma}{Lemma}

\newtheorem{conjecture}{Conjecture}
\theoremstyle{definition}

\theoremstyle{remark}
\newtheorem{remark}{Remark}

 \vspace{4mm}

\begin{document}

\hbox{Nanjing Univ. J. Math. Biquarterly 31(2014), no.\,2, 150--164.}
\medskip

\title
[{Some new series for $1/\pi$ and related congruences}]
{Some new series for $1/\pi$ and related congruences}

\author
[Zhi-Wei Sun] {Zhi-Wei Sun}

\address {Department of Mathematics, Nanjing
University, Nanjing 210093, People's Republic of China}
\email{zwsun@nju.edu.cn}

\keywords{Series for $1/\pi$, binomial coefficients, congruences.
\newline \indent The work was supported by the National Natural Science Foundation (grant 11171140) of China, and
the initial version of this paper was posted to {\tt arXiv} in April 2011 as a preprint with the ID {\tt arXiv:1104.3856}.}

 \begin{abstract} In this paper we prove some new series for $1/\pi$ as well as related congruences.
We also raise several new kinds of series for $1/\pi$ and
present some related conjectural congruences involving representations of primes by binary quadratic forms.
\end{abstract}

\maketitle

\section{Introduction}
\setcounter{lemma}{0}
\setcounter{theorem}{0}
\setcounter{corollary}{0}
\setcounter{remark}{0}
\setcounter{equation}{0}
\setcounter{conjecture}{0}

In 1914 S. Ramanujan \cite{R} (see also \cite{Be}) recorded 17 hypergeometric series for $1/\pi$ with six of them as follows:
\begin{align*}\sum_{k=0}^\infty\f{5k+1}{(-192)^k}\bi{2k}k^2\bi{3k}k=&\f 4{\sqrt3\,\pi},\tag{R1}
\\\sum_{k=0}^\infty\f{7k+1}{648^k}\bi{2k}k^2\bi{4k}{2k}=&\f 9{2\pi},\tag{R2}
\\\sum_{k=0}^\infty\f{65k+8}{(-63^2)^k}\bi{2k}k^2\bi{4k}{2k}=&\f {9\sqrt7}{\pi},\tag{R3}
\\\sum_{k=0}^\infty\f{6k+1}{256^k}\bi{2k}k^3=&\f4{\pi},\tag{R4}
\\\sum_{k=0}^\infty\f{6k+1}{(-512)^k}\bi{2k}k^3=&\f{2\sqrt2}{\pi},\tag{R5}
\\\sum_{k=0}^\infty\f{42k+5}{4096^k}\bi{2k}k^3=&\f{16}{\pi}.\tag{R6}
\end{align*}
For the techniques to prove them, see \cite{ChCh,BBC,BB}.
L. van Hamme \cite{vH} noted that such series usually have corresponding $p$-adic congruences.
In this paper we will present some new kinds of series for $1/\pi$ and give related congruences.

Our main results are as follows.

\begin{theorem}\label{Th1.1} {\rm (i)} We have
\begin{align} \label{1.1}\sum_{n=0}^\infty\f n{32^n}\sum_{k=0}^n\bi{2k}k^2\bi{2n-2k}{n-k}^2=&\f2{\pi},
\\\label{1.2}\sum_{n=0}^\infty\f n{54^n}\sum_{k=0}^n\bi{2k}k\bi{3k}k\bi{2(n-k)}{n-k}\bi{3(n-k)}{n-k}=&\f{\sqrt3}{\pi},
\\\label{1.3}\sum_{n=0}^\infty\f n{128^n}\sum_{k=0}^n\bi{4k}{2k}\bi{2k}k\bi{4(n-k)}{2(n-k)}\bi{2(n-k)}{n-k}=&\f{\sqrt2}{\pi},
\\\label{1.4}\sum_{n=0}^\infty\f n{864^n}\sum_{k=0}^n\bi{6k}{3k}\bi{3k}k\bi{6(n-k)}{3(n-k)}\bi{3(n-k)}{n-k}=&\f{1}{\pi}.
\end{align}

{\rm (ii)} Let $p$ be an odd prime. Then
\begin{align}\label{1.5}\sum_{n=0}^{p-1}\f n{32^n}\sum_{k=0}^n\bi{2k}k^2\bi{2(n-k)}{n-k}^2\eq&0\pmod{p^3},
\\\label{1.6}\sum_{n=0}^{p-1}\f n{128^n}\sum_{k=0}^n\bi{4k}{2k}\bi{2k}k\bi{4(n-k)}{2(n-k)}\bi{2(n-k)}{n-k}\eq&0\pmod{p^2}.
\end{align}
If $p>3$, then
\begin{align}\label{1.7}\sum_{n=0}^{p-1}\f n{54^n}\sum_{k=0}^n\bi{2k}{k}\bi{3k}k\bi{2(n-k)}{n-k}\bi{3(n-k)}{n-k}\eq&0\pmod{p^2},
\\\label{1.8}\sum_{n=0}^{p-1}\f n{864^n}\sum_{k=0}^n\bi{6k}{3k}\bi{3k}k\bi{6(n-k)}{3(n-k)}\bi{3(n-k)}{n-k}\eq&0\pmod{p^2}.
\end{align}
\end{theorem}
\begin{remark}\label{Rem1.1} For $n\in\N=\{0,1,2,\ldots\}$, define
$$a_n:=\f1{(2n-1)\bi{3n}n}\sum_{k=0}^n\bi{6k}{3k}\bi{3k}k\bi{6(n-k)}{3(n-k)}\bi{3(n-k)}{n-k}.$$
The author ever conjectured that $a_n\in\Z$ for all $n=0,1,2,\ldots$, and $\lim_{n\to+\infty}\root{n}\of{a_n}=64$.
This was confirmed by V.J.W. Guo \cite{G}.
\end{remark}

\begin{theorem}\label{Th1.2} We have
\begin{equation}\label{1.9}\sum_{n=0}^\infty\f{9n+2}{(-216)^n}\sum_{k=0}^n\bi{2k}k\bi{3k}k\bi{2(n-k)}{n-k}\bi{3(n-k)}{n-k}=\f 8{\sqrt3\,\pi}.
\end{equation}
Also,
\begin{equation}\label{1.10}\begin{aligned}&{}_2F_1\bigg(\begin{array}{c}1/3,2/3\\1\end{array}\bigg|{-\f18}\bigg)^2
\\=&\f14\,{}_2F_1\bigg(\begin{array}{c}1/3,2/3\\1\end{array}\bigg|{-\f18}\bigg)\,{}_2F_1\bigg(\begin{array}{c}4/3,5/3\\2\end{array}\bigg|{-\f18}\bigg)
+\f 4{\sqrt3\,\pi},
\end{aligned}\end{equation}
where
$${}_2F_1\bigg(\begin{array}{c}a,b\\c\end{array}\bigg|{z}\bigg):=1+\sum_{n=1}^\infty\(\prod_{k=0}^{n-1}\f{(a+k)(b+k)}{(k+1)(c+k)}\)z^n.$$
\end{theorem}

\begin{theorem}\label{Th1.3} We have
\begin{align}\label{1.11}\sum_{n=0}^\infty\f n{128^n}\bi{2n}np_n(4)=&\f{\sqrt2}{\pi},
\\\label{1.12}\sum_{n=0}^\infty\f{8n+1}{576^n}\bi{2n}np_n(4)=&\f 9{2\pi},
\\\label{1.13}\sum_{n=0}^\infty\f{8n+1}{(-4032)^n}\bi{2n}np_n(4)=&\f{9\sqrt7}{8\pi},
\end{align}
where
\begin{equation}\label{1.14}p_n(x):=\sum_{k=0}^n\bi{2k}k^2\bi{2(n-k)}{n-k}x^{n-k}\ \ (n=0,1,2,\ldots).\end{equation}
\end{theorem}
\begin{remark}\label{Rem1.2} In Section 4 we will present many conjectural series for $1/\pi$ involving $p_n(x)$.
\end{remark}
\medskip

\begin{theorem}\label{Th1.4} We have
\begin{align}\label{1.15} \sum_{n=0}^\infty\f n{4^n}\sum_{k=0}^n\bi{-1/4}k^2\bi{-3/4}{n-k}^2=&\f{4\sqrt3}{9\pi},
\\\label{1.16}\sum_{n=0}^\infty\f{9n+2}{(-8)^n}\sum_{k=0}^n\bi{-1/4}k^2\bi{-3/4}{n-k}^2=&\f 4{\pi},
\\\label{1.17}\sum_{n=0}^\infty\f{9n+1}{64^n}\sum_{k=0}^n\bi{-1/4}k^2\bi{-3/4}{n-k}^2=&\f {64}{7\sqrt7\,\pi}.
\end{align}
\end{theorem}

There are many series for $1/\pi$ of the form
$$\sum_{k=0}^\infty (bk+c)\f{\bi{2k}ka_k}{m^k}=\f {C}{\pi},$$
where $a_k,b,c,C$ and $m\not=0$ are real numbers. We note that if $m>8$ or $m<0$ then
\begin{equation}\label{1.18}\sum_{n=0}^\infty(bmn+2b+(m-4)c)\f{\bi{2n}na_n^*}{(4-m)^n}
=(m-4)\sqrt{\f{m-4}m}\sum_{k=0}^\infty(bk+c)\f{\bi{2k}ka_k}{m^k},
\end{equation}
where $a_n^*=\sum_{k=0}^n\bi nk(-1)^ka_k$.
This is easy, because
\begin{align*} &\sum_{n=0}^\infty(bmn+2b+(m-4)c)\f{\bi{2n}na_n^*}{(4-m)^n}
\\=&\sum_{n=0}^\infty(bmn+2b+(m-4)c)\f{\bi{2n}n}{(4-m)^n}\sum_{k=0}^n\bi nk(-1)^ka_k
\\=&\sum_{k=0}^\infty(-1)^ka_k\sum_{n=k}^\infty(bmn+2b+(m-4)c)\f{\bi{2n}n\bi nk}{(4-m)^n}
\\=&\sum_{k=0}^\infty(-1)^ka_k\((m-4)\sqrt{\f{m-4}m}\ (bk+c)\f{\bi{2k}k}{(-m)^k}\)
\\&\l(\t{since}\ \bi{2n}n\bi nk=(-4)^n\bi{-1/2}n\bi nk=(-4)^{n-k}\bi{2k}k\bi{-1/2-k}{n-k}\r)
\\=&(m-4)\sqrt{\f{m-4}m}\sum_{k=0}^\infty(bk+c)\f{\bi{2k}ka_k}{m^k}.
\end{align*}

We will show Theorems 1.1-1.3 in the next section, and prove Theorem 1.4 in Section 3.
Section 4 contains some new kinds of conjectures on series for $1/\pi$ or related congruences.
Note that we have explained in \cite{Su4} and \cite{Su5} how the author found some new series for $1/\pi$.

\section{Proofs of Theorems 1.1-1.3}
\setcounter{lemma}{0}
\setcounter{theorem}{0}
\setcounter{corollary}{0}
\setcounter{remark}{0}
\setcounter{equation}{0}
\setcounter{conjecture}{0}

\begin{lemma}\label{Lem2.1} {\rm (\cite[Lemma 3.1]{Su3})} For any $n=0,1,2,\ldots$ we have
\begin{equation}\label{2.1}\sum_{k=0}^n\bi{2k}k^3\bi k{n-k}(-16)^{n-k}=\sum_{k=0}^n\bi{2k}k^2\bi{2(n-k)}{n-k}^2.
\end{equation}
\end{lemma}
\begin{remark}\label{Rem2.1} Z.-H. Sun \cite{S} proved that for any $n\in\N$ we have
\begin{equation}\label{2.2}\begin{aligned}&\sum_{k=0}^n\bi{2k}{k}^2\bi{3k}k\bi k{n-k}(-27)^{n-k}
\\=&\sum_{k=0}^n\bi{2k}{k}\bi{3k}k\bi{2(n-k)}{n-k}
\bi{3(n-k)}{n-k},\end{aligned}\end{equation}
\begin{equation}\label{2.3}\begin{aligned}&\sum_{k=0}^n\bi{4k}{2k}\bi{2k}k^2\bi k{n-k}(-64)^{n-k}
\\=&\sum_{k=0}^n\bi{4k}{2k}\bi{2k}k\bi{4(n-k)}{2(n-k)}
\bi{2(n-k)}{n-k},
\end{aligned}\end{equation}
\begin{equation}\label{2.4}\begin{aligned}&\sum_{k=0}^n\bi{6k}{3k}\bi{3k}k\bi{2k}k\bi k{n-k}(-432)^{n-k}
\\=&\sum_{k=0}^n\bi{6k}{3k}\bi{3k}k\bi{6(n-k)}{3(n-k)}
\bi{3(n-k)}{n-k}.
\end{aligned}\end{equation}
Actually (\ref{2.1})-(\ref{2.4}) can be easily deduced from the classical Clausen identity (cf. \cite[p.\,116]{AAR}) for hypergeometric series.
\end{remark}
\medskip
\noindent{\it Proof of Theorem 1.1}. Let $\sigma$ denote the left-hand side of (\ref{1.1}). Then
\begin{align*} \sigma=&\sum_{k=0}^\infty\bi{2k}k^2\sum_{n=k}^\infty\f n{32^n}\bi{2(n-k)}{n-k}^2
\\=&\sum_{k=0}^\infty\f{\bi{2k}k^2}{32^k}\sum_{j=0}^\infty\f{k+j}{32^j}\bi{2j}j^2.
\end{align*}
By {\tt Mathematica},
$$\sum_{k=0}^\infty\f{\bi{2k}k^2}{32^k}=\f{\sqrt{\pi}}{\Gamma(3/4)^2}
\ \ \t{and}\ \
\sum_{k=0}^\infty\f{k\bi{2k}k^2}{32^k}=\f{2\sqrt{\pi}}{\Gamma(1/4)^2},$$
where $\Gamma(x)$ is the well-known $\Gamma$-function.
By a known forumla,
$$\Gamma\l(\f14\r)\Gamma\l(\f 34\r)=\f{\pi}{\sin(\pi/4)}=\sqrt2\,\pi.$$
So we have
\begin{align*} \sigma=&\sum_{k=0}^\infty\f{\bi{2k}k^2}{32^k}\(k\f{\sqrt{\pi}}{\Gamma(3/4)^2}+\f{2\sqrt{\pi}}{\Gamma(1/4)^2}\)
\\=&\f{\sqrt{\pi}}{\Gamma(3/4)^2}\cdot\f{2\sqrt{\pi}}{\Gamma(1/4)^2}
+\f{2\sqrt{\pi}}{\Gamma(1/4)^2}\cdot\f{\sqrt{\pi}}{\Gamma(3/4)^2}
\\=&\f{2\pi}{2\pi^2}+\f{2\pi}{2\pi^2}=\f2{\pi}.
\end{align*}

Similarly, using {\tt Mathematica} we find
\begin{gather*}\sum_{k=0}^\infty\f{\bi{2k}k\bi{3k}k}{54^k}=\f{\sqrt{\pi}}{\Gamma(2/3)\Gamma(5/6)}
\ \ \t{and}\ \
\sum_{k=0}^\infty\f{k\bi{2k}k\bi{3k}k}{54^k}=\f{2\sqrt{\pi}}{\Gamma(1/6)\Gamma(1/3)},
\\\sum_{k=0}^\infty\f{\bi{4k}{2k}\bi{2k}k}{128^k}=\f{\sqrt{\pi}}{\Gamma(5/8)\Gamma(7/8)}
\ \ \t{and}\ \
\sum_{k=0}^\infty\f{k\bi{4k}{2k}\bi{2k}k}{128^k}=\f{2\sqrt{\pi}}{\Gamma(1/8)\Gamma(3/8)},
\\\sum_{k=0}^\infty\f{\bi{6k}{3k}\bi{3k}k}{864^k}=\f{\sqrt{\pi}}{\Gamma(7/12)\Gamma(11/12)}
\ \ \t{and}\ \
\sum_{k=0}^\infty\f{k\bi{6k}{3k}\bi{3k}k}{864^k}=\f{2\sqrt{\pi}}{\Gamma(1/12)\Gamma(5/12)}.
\end{gather*}
Thus we can prove (\ref{1.2})-(\ref{1.4}) in the way we show (\ref{1.1}).

(ii) Let $S$ denote the left-hand side of (\ref{1.5}). By Lemma \ref{Lem2.1}, we have
\begin{align*} S=&\sum_{n=0}^{p-1}\f n{32^n}\sum_{k=0}^n\bi{2k}k^3\bi k{n-k}(-16)^{n-k}
\\=&\sum_{k=0}^{p-1}\f{\bi{2k}k^3}{32^k}\sum_{j=0}^{p-1-k}(k+j)\bi kj\f{(-16)^j}{32^j}.
\end{align*}
If $(p-1)/2<k\ls p-1$, then $p\mid \bi{2k}k$ and hence $\bi{2k}k^3\eq0\pmod{p^3}$.
If $0\ls k\ls(p-1)/2$, then $p-1-k\gs k$ and hence
\begin{align*}&\sum_{j=0}^{p-1-k}(k+j)\bi kj\f1{(-2)^j}
\\=&k\sum_{j=0}^k\bi kj\f 1{(-2)^j}+\f k{-2}\sum_{j=1}^{k-1}\bi{k-1}{j-1}\f1{(-2)^{j-1}}
\\=&\f k{2^k}-\f k2\cdot\f1{2^{k-1}}=0.
\end{align*}
So, by the above we have $S\eq0\pmod{p^3}$.

Similarly, by using (\ref{2.2})-(\ref{2.4}) we can easily prove (\ref{1.6})-(\ref{1.8}).

The proof of Theorem 1.1 is now complete. \qed

\medskip
\noindent{\it Proof of Theorem 1.2}. Let $S$ denote the left-hand side of (\ref{1.9}). In view of (\ref{2.2}), we have
\begin{align*} S=&\sum_{n=0}^\infty\f{9n+2}{(-216)^n}\sum_{k=0}^n\bi{2k}k^2\bi{3k}k\bi k{n-k}(-27)^{n-k}
\\=&\sum_{k=0}^\infty\f{\bi{2k}k^2\bi{3k}k}{(-216)^k}\sum_{j=0}^k(9(k+j)+2)\bi kj\f{(-27)^j}{(-216)^j}
\\=&\sum_{k=0}^\infty\f{\bi{2k}k^2\bi{3k}k}{(-216)^k}\((9k+2)\sum_{j=0}^k\bi kj\f1{8^j}+\f{9k}8\sum_{j=1}^{k-1}\bi{k-1}{j-1}\f1{8^{j-1}}\)
\\=&\sum_{k=0}^\infty\f{\bi{2k}k^2\bi{3k}k}{(-216)^k}(9k+2+k)\f{9^k}{8^k}=2\sum_{k=0}^\infty\f{5k+1}{(-192)^k}\bi{2k}k^2\bi{3k}k
\\=&2\times\f 4{\sqrt3\,\pi}=\f 8{\sqrt3\,\pi}\ \ \ (\t{by Ramanujan's series (R1)}).
\end{align*}
This proves (\ref{1.9}). On the other hand,
\begin{align*} S=&\sum_{k=0}^\infty\f{\bi{2k}k\bi{3k}k}{(-216)^k}\sum_{n=k}^\infty\f{9n+2}{(-216)^{n-k}}\bi{2(n-k)}{n-k}\bi{3(n-k)}{n-k}
\\=&\sum_{k=0}^\infty\f{\bi{2k}k\bi{3k}k}{(-216)^k}\sum_{j=0}^\infty\f{9(k+j)+2}{(-216)^j}\bi{2j}j\bi{3j}j
\\=&\sum_{k=0}^\infty\f{\bi{2k}k\bi{3k}k}{(-216)^k}\((9k+2)\sum_{j=0}^\infty\f{\bi{2j}j\bi{3j}j}{(-216)^j}+9\sum_{j=0}^\infty\f{j\bi{2j}j\bi{3j}j}{(-216)^j}\)
\\=&18\sum_{k=0}^\infty\f{k\bi{2k}k\bi{3k}k}{(-216)^k}\sum_{j=0}^\infty\f{\bi{2j}j\bi{3j}j}{(-216)^j}
+2\(\sum_{k=0}^\infty\f{\bi{2k}k\bi{3k}k}{(-216)^k}\)^2.
\end{align*}
Note that
$$\sum_{k=0}^\infty\f{\bi{2k}k\bi{3k}k}{(-216)^k}={}_2F_1\bigg(\begin{array}{c}1/3,2/3\\1\end{array}\bigg|{-\f18}\bigg)$$
and
$$\sum_{k=0}^\infty\f{k\bi{2k}k\bi{3k}k}{(-216)^k}=-\f1{36}\,{}_2F_1\bigg(\begin{array}{c}4/3,5/3\\2\end{array}\bigg|{-\f18}\bigg).$$
So we also have (\ref{1.10}). This concludes the proof. \qed

\begin{lemma}\label{Lem2.2} For all $n\in\N$ we have
\begin{equation}\label{2.5}\bi{2n}np_n(4)
=\sum_{k=0}^n\bi{4k}{2k}\bi{2k}k^2\bi k{n-k}(-64)^{n-k}.
\end{equation}
\end{lemma}
\Proof. Define
$$u_n=\f1{\bi{2n}n}\sum_{k=0}^n\bi{4k}{2k}\bi{2k}k^2\bi k{n-k}(-64)^{n-k}$$
for $n=0,1,2,\ldots$. It is easy to see that $u_0=p_0(4)=1$ and $u_1=p_1(4)=12$.
By the Zeilberger algorithm (cf. \cite[pp.\,101-119]{PWZ}), we find via {\tt Mathematica} that
$$(n+2)^2u_{n+2}=4(8n^2+24n+19)u_{n+1}-256(n+1)^2u_n$$
and
$$(n+2)^2p_{n+2}(4)=4(8n^2+24n+19)p_{n+1}(4)-256(n+1)^2p_n(4)$$
for all $n=0,1,2,\ldots$. Therefore, $u_n=p_n(4)$ for all $n\in\N$.
\qed

\medskip
\noindent{\it Proof of Theorem 1.3}.  By (\ref{2.3}) and (\ref{2.5}), (\ref{1.11}) is equivalent to (\ref{1.3}).

In view of (\ref{2.4}) and (\ref{2.5}),
\begin{align*}&\sum_{n=0}^\infty\f{8n+1}{576^n}\bi{2n}np_n(4)
\\=&\sum_{n=0}^\infty\f{8n+1}{576^n}\sum_{k=0}^n\bi{4k}{2k}\bi{2k}k^2\bi k{n-k}(-64)^{n-k}
\\=&\sum_{k=0}^\infty\f{\bi{4k}{2k}\bi{2k}k^2}{576^k}\sum_{j=0}^k(8(k+j)+1)\bi kj\l(\f{-64}{576}\r)^j
\\=&\sum_{k=0}^\infty\f{\bi{4k}{2k}\bi{2k}k^2}{576^k}\((8k+1)\sum_{j=0}^k\bi kj\f1{(-9)^j}-\f {8k}9\sum_{j=1}^k\bi{k-1}{j-1}\f1{(-9)^{j-1}}\)
\\=&\sum_{k=0}^\infty\f{\bi{4k}{2k}\bi{2k}k^2}{576^k}((8k+1)-k)\l(1-\f19\r)^k
\\=&\sum_{k=0}^\infty(7k+1)\f{\bi{4k}{2k}\bi{2k}k^2}{648^k}=\f 9{2\pi}\ \ \t{(by (R2))}.
\end{align*}
So (\ref{1.12}) holds. Similarly, by using (R3) we can prove (\ref{1.13}).

The proof of Theorem 1.3 is now complete. \qed

\section{Proof of Theorem 1.4}
\setcounter{lemma}{0}
\setcounter{theorem}{0}
\setcounter{corollary}{0}
\setcounter{remark}{0}
\setcounter{equation}{0}
\setcounter{conjecture}{0}

\begin{lemma}\label{Lem3.1} For any $n\in\N$ we have
\begin{equation}\label{3.1}64^n\sum_{k=0}^n\bi{-1/4}k^2\bi{-3/4}{n-k}^2=\sum_{k=0}^n\bi{2k}k^3\bi{2(n-k)}{n-k}16^{n-k}.\end{equation}
\end{lemma}
\Proof. Let $u_n$ denote the left-hand side or the right-hand side of (\ref{3.1}). It is easy to see that $u_0=1$ and $u_1=40$.
By the Zeilberger algorithm we find the recurrence relation
$$(n+2)^3u_{n+2}=8(2n+3)(8n^2+24n+21)u_{n+1}-4096(n+1)^3u_n\ (n\in\N).$$
So, by induction, (\ref{3.1}) holds for all $n=0,1,2,\ldots$. \qed

\begin{lemma}\label{Lem3.2} For any number $m$ with $|m|>4$, we have
\begin{equation}\label{3.2}\sum_{k=0}^\infty\f{\bi{2k}k}{m^k}=\sqrt{\f m{m-4}}
\ \ \t{and}\ \ \sum_{k=0}^\infty\f{k\bi{2k}k}{m^k}=\f2{m-4}\sqrt{\f m{m-4}}.
\end{equation}
\end{lemma}
\Proof. Clearly $\bi{2k}k=\bi{-1/2}k(-4)^k$ for all $k\in\N$. Thus
$$\sum_{k=0}^\infty\f{\bi{2k}k}{m^k}=\sum_{k=0}^\infty\bi{-1/2}k\l(-\f 4m\r)^k=\l(1-\f 4m\r)^{-1/2}=\sqrt{\f m{m-4}}$$
and
\begin{align*} \sum_{k=0}^\infty\f{k\bi{2k}k}{m^k}=&\sum_{k=1}^\infty k\bi{-1/2}k\l(-\f 4m\r)^k
\\=&-\f12\l(-\f 4m\r)\sum_{k=1}^\infty\bi{-1/2-1}{k-1}\l(-\f 4m\r)^{k-1}
\\=&\f2m\l(1-\f 4m\r)^{-3/2}=\f2{m-4}\sqrt{\f m{m-4}}.
\end{align*}
This concludes the proof. \qed

\medskip
\noindent{\it Proof of Theorem 1.4}.
We reduce (\ref{1.15})-(\ref{1.17}) to (R4)-(R6) respectively.
Below we just give a detailed proof of (\ref{1.17}). (\ref{1.15}) and (\ref{1.16}) can be proved in a similar way.

Let $S$ denote the left-hand side of (\ref{1.17}). In view of Lemmas 3.1-3.2 and (R6),
\begin{align*} S=&\sum_{n=0}^\infty\f{9n+1}{64^{2n}}\sum_{k=0}^n\bi{2k}k^3\bi{2(n-k)}{n-k}16^{n-k}
\\=&\sum_{k=0}^\infty\f{\bi{2k}k^3}{4096^k}\sum_{j=0}^\infty\f{9(k+j)+1}{256^j}\bi{2j}j
\\=&\sum_{k=0}^\infty\f{\bi{2k}k^3}{4096^k}\l((9k+1)\sqrt{\f{256}{252}}+9\times\f2{252}\sqrt{\f{256}{252}}\r)
\\=&\f4{7\sqrt7}\sum_{k=0}^\infty\f{42k+5}{4096^k}\bi{2k}k^3=\f4{7\sqrt7}\times\f{16}{\pi}=\f{64}{7\sqrt7\,\pi}.
\end{align*}
So (\ref{1.17}) holds. \qed

\section{Some related conjectures}
\setcounter{lemma}{0}
\setcounter{theorem}{0}
\setcounter{corollary}{0}
\setcounter{remark}{0}
\setcounter{equation}{0}
\setcounter{conjecture}{0}

Recall that the Euler numbers $E_0,E_1,E_2,\ldots$ are integers defined by
$$E_0=1\ \ \t{and}\ \ \sum^n_{k=0\atop2\mid k}\bi nk E_{n-k}=0\ \ (n=1,2,3,\ldots).$$

In view of Theorem 1.1, we raise the following conjecture.

\begin{conjecture}\label{Conj4.1} Let $p$ be an odd prime. Then
$$\sum_{n=0}^{p-1}\f n{32^n}\sum^n_{k=0}\bi{2k}k^2\bi{2n-2k}{n-k}^2\eq-2p^3E_{p-3}\pmod{p^4}.$$
When $p>3$, we also have
\begin{align*}\sum_{k=0}^{p-1}\f{n+1}{8^n}\sum_{k=0}^n\bi{2k}k^2\bi{2n-2k}{n-k}^2\eq& (-1)^{(p-1)/2}p+5p^3E_{p-3}\pmod{p^4},
\\\sum_{k=0}^{p-1}\f{2n+1}{(-16)^n}\sum_{k=0}^n\bi{2k}k^2\bi{2n-2k}{n-k}^2\eq& (-1)^{(p-1)/2}p+3p^3E_{p-3}\pmod{p^4}.
\end{align*}
\end{conjecture}
\begin{remark}\label{Rem4.1}. In \cite{Su1,Su2}, the author established some basic $p$-adic congruences involving $E_{p-3}$.
\end{remark}
\medskip

Motivated by Theorem 1.3, we pose the following conjecture.

\begin{conjecture}\label{Conj4.2} We have
\begin{align}\label{4.1}\sum_{k=0}^\infty\f{k-1}{72^k}\bi{2k}k p_k(4)=&\f{9}{\pi},
\\\label{4.2}\sum_{k=0}^\infty\f{4k+1}{(-192)^k}\bi{2k}k p_k(4)=&\f{\sqrt{3}}{\pi},
\\\label{4.3}\sum_{k=0}^{\infty}\f{k-2}{100^k}\bi{2k}k p_k(6)=&\f{50}{3\pi},
\\\label{4.4}\sum_{k=0}^\infty\f{k}{(-192)^k}\bi{2k}k p_k(-8)=&\f{3}{2\pi},
\\\label{4.5}\sum_{k=0}^\infty\f{6k-1}{256^k}\bi{2k}kp_k(12)=&\f{8\sqrt3}{\pi},
\\\label{4.6}\sum_{k=0}^\infty\f{17k-224}{(-225)^k}\bi{2k}kp_k(-14)=&\f{1800}{\pi},
\\\label{4.7}\sum_{k=0}^\infty\f{15k-256}{289^k}\bi{2k}kp_k(18)=&\f{2312}{\pi},
\\\label{4.8}\sum_{k=0}^\infty\f{20k-11}{(-576)^k}\bi{2k}kp_k(-32)=&\f{90}{\pi},
\\\label{4.9}\sum_{k=0}^\infty\f{10k+1}{(-1536)^k}\bi{2k}kp_k(-32)=&\f{3\sqrt{6}}{\pi},
\\\label{4.10}\sum_{k=0}^\infty\f{3k-2}{640^k}\bi{2k}kp_k(36)=&\f{5\sqrt{10}}{\pi},
\\\label{4.11}\sum_{k=0}^\infty\f{12k+1}{1600^k}\bi{2k}kp_k(36)=&\f{75}{8\pi},
\\\label{4.12}\sum_{k=0}^\infty\f{24k+5}{3136^k}\bi{2k}kp_k(-60)=&\f{49\sqrt{3}}{8\pi},
\\\label{4.13}\sum_{k=0}^\infty\f{14k+3}{(-3072)^k}\bi{2k}kp_k(64)=&\f{6}{\pi},
\\\label{4.14}\sum_{k=0}^\infty\f{20k-67}{(-3136)^k}\bi{2k}kp_k(-192)=&\f{490}{\pi},
\end{align}
\begin{align}
\label{4.15}\sum_{k=0}^\infty\f{7k-24}{3200^k}\bi{2k}kp_k(196)=&\f{125\sqrt2}{\pi},
\\\label{4.16}\sum_{k=0}^\infty\f{5k-32}{(-6336)^k}\bi{2k}kp_k(-392)=&\f{495}{2\pi},
\\\label{4.17}\sum_{k=0}^\infty\f{66k-427}{6400^k}\bi{2k}kp_k(396)=&\f{1000\sqrt{11}}{\pi},
\\\label{4.18}\sum_{k=0}^\infty\f{34k-7}{(-18432)^k}\bi{2k}kp_k(-896)=&\f{54\sqrt{2}}{\pi},
\\\label{4.19}\sum_{k=0}^\infty\f{24k-5}{18496^k}\bi{2k}kp_k(900)=&\f{867}{16\pi}.
\end{align}
\end{conjecture}
\begin{remark}\label{Rem4.2} (\ref{4.5}) also appeared in \cite[Conjecture 1.5]{Su4}.
\end{remark}
\begin{conjecture}\label{Conj4.3} For $n=0,1,2,\ldots$ define
$$S_n(x)=\sum_{k=0}^n\bi nk\bi{2k}k\bi{2n-2k}{n-k}x^{n-k}.$$
Then we have
\begin{align}\label{4.20}\sum_{k=0}^\infty\f{12k+1}{400^k}\bi{2k}kS_k(16)=&\f{25}{\pi},
\\\label{4.21}\sum_{k=0}^\infty\f{10k+1}{(-384)^k}\bi{2k}kS_k(-16)=&\f{8\sqrt6}{\pi},
\\\label{4.22}\sum_{k=0}^\infty\f{170k+37}{(-3584)^k}\bi{2k}kS_k(64)=&\f{64\sqrt{14}}{3\pi},
\\\label{4.23}\sum_{k=0}^\infty\f{476k+103}{3600^k}\bi{2k}kS_k(-64)=&\f{225}{\pi},
\\\label{4.24}\sum_{k=0}^\infty\f{140k+19}{4624^k}\bi{2k}kS_k(64)=&\f{289}{3\pi},
\\\label{4.25}\sum_{k=0}^\infty\f{1190k+163}{(-4608)^k}\bi{2k}kS_k(-64)=&\f{576\sqrt2}{\pi}.
\end{align}
\end{conjecture}
\begin{remark}\label{Rem4.3} Note that $$S_n(-1)=\begin{cases}\bi{n}{n/2}^2&\t{if}\ 2\mid n,\\0&\t{if}\ 2\nmid n.\end{cases}$$
Also,
$$S_n(1)=\sum_{k=0}^{\lfloor n/2\rfloor}\bi n{2k}\bi{2k}k^24^{n-2k}.$$
The two identities can be easily proved via the Zeilberger algorithm.
Identities of the form
$$\sum_{n=0}^\infty \f{bn+c}{m^n}\bi{2n}n\sum_{k=0}^{\lfloor n/2\rfloor}\bi n{2k}\bi{2k}k^24^{n-2k}=\f C{\pi}$$
were recently investigated in \cite{CC}.
\end{remark}

\begin{conjecture}\label{Conj4.4} Define
$$s_n(x):=\sum_{k=0}^n\bi nk\bi{n+2k}{2k}\bi{2k}kx^{-(n+k)}\quad\t{for}\ n=0,1,2,\ldots$$
Then
\begin{align}\label{4.26}\sum_{k=0}^\infty(7k+2)\bi{2k}ks_k(-9)&=\f{9\sqrt3}{5\pi},
\\\label{4.27}\sum_{k=0}^\infty(9k+2)\bi{2k}ks_k(-20)&=\f{4}{\pi},
\\\label{4.28}\sum_{k=0}^\infty(95k+13)\bi{2k}ks_k(36)&=\f{18\sqrt{15}}{\pi},
\\\label{4.29}\sum_{k=0}^\infty(310k+49)\bi{2k}ks_k(-64)&=\f{32\sqrt{15}}{\pi},
\\\label{4.30}\sum_{k=0}^\infty(495k+53)\bi{2k}ks_k(196)&=\f{70\sqrt{7}}{\pi},
\\\label{4.31}\sum_{k=0}^\infty(13685k+1474)\bi{2k}ks_k(-324)&=\f{1944\sqrt5}{\pi},
\\\label{4.32}\sum_{k=0}^\infty(3245k+268)\bi{2k}ks_k(1296)&=\f{1215}{\sqrt2\,\pi},
\\\label{4.33}\sum_{k=0}^\infty(6420k+443)\bi{2k}ks_k(5776)&=\f{1292\sqrt{95}}{9\pi}.
\end{align}
Also,
\begin{align}\label{4.34}\sum_{n=0}^\infty\f{357n+103}{2160^n}\bi{2n}n\sum_{k=0}^n\bi nk\bi{n+2k}{2k}\bi{2k}k(-324)^{n-k}=&\f{90}{\pi},
\\\label{4.35}\sum_{n=0}^\infty\f{n}{3645^n}\bi{2n}n\sum_{k=0}^n\bi nk\bi{n+2k}{2k}\bi{2k}k486^{n-k}=&\f{10}{3\pi}.
\end{align}
\end{conjecture}
\begin{remark}\label{Rem4.4} (4.34) also appeared in \cite[Conjecture 1.7]{Su4}. The use of the polynomials $\sum_{k=0}^n\bi nk\bi{n+2k}{2k}\bi{2k}kx^{n-k}\ (n\in\N)$
 was inspired by an identity of MacMahon \cite[p.\,122]{M} which states that
$$\sum_{k=0}^n\bi nk\bi{n+2k}{2k}\bi{2k}k(-4)^{n-k}=\sum_{k=0}^n\bi nk^3\quad\t{for all}\ n=0,1,2,\ldots.$$
\end{remark}
\medskip

Conjectures 4.2-4.4 were motivated by our investigation of related
congruences. For example, (\ref{4.35}) was
inspired by the following conjecture.

\begin{conjecture}\label{Conj4.5} Let $p>5$ be a prime. Then
\begin{align*}&\sum_{n=0}^{p-1}\f{\bi{2n}n}{3645^n}\sum_{k=0}^n\bi nk\bi{n+2k}{2k}\bi{2k}k486^{n-k}
\\\eq&\begin{cases} 4x^2-2p\pmod{p^2}&\t{if}\ p\eq1,4\pmod{15}\ \&\ p=x^2+15y^2\ (x,y\in\Z),
\\2p-12x^2\pmod{p^2}&\t{if}\ p\eq2,8\pmod{15}\ \&\ p=3x^2+5y^2\ (x,y\in\Z),
\\0\pmod{p^2}&\t{if}\ (\f p{15})=-1,\end{cases}
\end{align*}
where $(-)$ denotes the Jacobi symbol.
Also,
\begin{align*}&\sum_{n=0}^{p-1}\f{n\bi{2n}n}{3645^n}\sum_{k=0}^n\bi nk\bi{n+2k}{2k}\bi{2k}k486^{n-k}
\\&\qquad\eq\f{16p}{81}\l(\l(\f{-1}p\r)-\l(\f 5p\r)\r)\pmod{p^2}.
\end{align*}
\end{conjecture}
\begin{remark}\label{Rem4.5} For any prime $p>5$ with $(\f p{15})=1$, it is known (see, e.g., \cite{Co}) that $p$ can be written
in the form $x^2+15y^2$ with $x,y\in\Z$ (or the form $3x^2+5y^2\ (x,y\in\Z))$ if $p\eq1,4\pmod{15}$ (or $p\eq2,8\pmod{15}$, respectively).
\end{remark}

Now we pose a very curious conjecture.

\begin{conjecture}\label{Conj4.6} Let $p>3$ be a prime, and set
$$R_p:=\f1p\sum_{n=0}^{p-1}\f{6n+1}{(-1728)^n}\bi{2n}n\sum_{k=0}^n\bi nk\bi{n+2k}{2k}\bi{2k}k(-324)^{n-k}.$$
Then
$$R_p^2\eq\f{512(\f{10}p)-27(\f{-15}p)-475}{10}\ \pmod p.$$
If $p>5$, then
\begin{align*}&\sum_{n=0}^{p-1}\f{\bi{2n}n}{(-1728)^n}\sum_{k=0}^n\bi nk\bi{n+2k}{2k}\bi{2k}k(-324)^{n-k}
\\\eq&\begin{cases}4x^2-2p\pmod{p^2}&\t{if}\ p\eq1,7\pmod{24}\ \&\ p=x^2+6y^2\,(x,y\in\Z),
\\8x^2-2p\pmod{p^2}&\t{if}\ p\eq5,11\pmod{24}\ \&\ p=2x^2+3y^2\,(x,y\in\Z),
\\0\pmod{p^2}&\t{if}\ (\f{-6}p)=-1,\ \t{i.e.},\ p\eq13,17,19,23\pmod{24}.
\end{cases}
\end{align*}
\end{conjecture}
\begin{remark}\label{Rem4.6} For any prime $p>3$ with $(\f {-6}p)=1$, it is known (cf. \cite{Co}) that $p$ can be written
in the form $x^2+6y^2$ with $x,y\in\Z$ (or the form $2x^2+3y^2\ (x,y\in\Z))$ if $p\eq1,7\pmod{24}$ (or $p\eq5,11\pmod{24}$, respectively).
\end{remark}

A sequence of polynomials $\{P_n(q)\}_{n\gs0}$ with integer coefficients is said to be {\it $q$-logconvex} if
for each  $n=1,2,3,\ldots$ all the coefficients of the polynomial $P_{n-1}(q)P_{n+1}(q)-P_n(q)^2\in\Z[q]$ are nonnegative.
In view of Conjecture 4.3 and \cite[Section 6]{Su4}, we propose the following conjecture.

\begin{conjecture}\label{Conj4.7} $\{P_n(q)\}_{n\gs0}$ is $q$-logconvex if $P_n(q)$ has one of the following three forms:
$$\sum_{k=0}^n\bi nk^2\bi{n+k}kq^k,\ \ \ S_n(q)=\sum_{k=0}^n\bi nk\bi{2k}k\bi{2(n-k)}{n-k}q^k,$$
and $$D_n(q)=\sum_{k=0}^n\bi nk^2\bi{2k}k\bi{2(n-k)}{n-k}q^k.$$
\end{conjecture}


\begin{thebibliography}{99}

\bibitem{AAR} G. E. Andrews, R. Askey and R. Roy, Special Functions, Cambridge Univ. Press, Cambridge, 1999.
\bibitem{BB} N. D. Baruah and B. C. Berndt, {\it Eisenstein series and Ramanujan-type series for $1/\pi$},
 Ramanujan J. {\bf 23} (2010), 17--44.

\bibitem{BBC} N. D. Baruah, B. C. Berndt and H. H. Chan,
{\it Ramanujan' series for $1/\pi$: a survey}, Amer. Math. Monthly {\bf 116} (2009), 567--587.

\bibitem{Be} B. C. Berndt,  Ramanujan's Notebooks, Part IV, Springer, New York, 1994.

\bibitem{CC} H. H. Chan and S. Cooper, {\it Rational analogues of Ramanujan's series for $1/\pi$},
Math. Proc. Cambridge Philos. Soc. {\bf 153} (2012), 361--383.

\bibitem{ChCh} D. V. Chudnovsky and G. V. Chudnovsky, {\it Approximations and complex multiplication
according to Ramanujan}, in: Ramanujan Revisited: Proc. of the Centenary Confer. (Urbana-Champaign, ILL., 1987),
(eds., G. E. Andrews, B. C. Berndt and R. A. Rankin), Academic Press, Boston, MA, 1988, pp. 375--472.

\bibitem{C} S. Cooper, {\it Sporadic sequences, modular forms and new series for $1/\pi$},
Ramanujan J. {\bf 29} (2012), 163--183.

\bibitem{Co} D. A. Cox, Primes of the Form $x^2+ny^2$, John Wiley \& Sons, 1989.

\bibitem{G} V.J.W. Guo, {\it Proof of Sun's conjecture on the divisibility of certain binomial sums},
Electron. J. Combin. {\bf 20} (2013), no. 4, \#P20.

\bibitem{M} P. A. MacMachon, Combinatorial Analysis, Vol. 1, Cambridge Univ. Press, London, 1915.

\bibitem{PWZ} M. Petkov\v sek, H. S. Wilf and D. Zeilberger, $A=B$, A K Peters, Wellesley, 1996.

\bibitem{R} S. Ramanujan, {\it Modular equations and approximations to $\pi$},
Quart. J. Math. (Oxford) {\bf (2)45} (1914), 350--372.

\bibitem{S} Z.-H. Sun, {\it Some supercongruences modulo $p^2$}, preprint, {\tt arXiv:1101.1050}.

\bibitem{Su1} Z.-W. Sun, {\it On congruences related to central binomial coefficients},
J. Number Theory {\bf 131} (2011), 2219--2238.

\bibitem{Su2} Z.-W. Sun, {\it Super congruences and Euler numbers},
Sci. China Math. {\bf 54} (2011), 2509--2535.

\bibitem{Su3} Z.-W. Sun, {\it On sums involving products of three binomial coefficients},
Acta Arith. {\bf 156} (2012), 123-141.

\bibitem{Su4} Z.-W. Sun, {\it Conjectures and results on $x^2$ mod $p^2$ with $4p=x^2+dy^2$},  in: Y. Ouyang, C. Xing, F.
Xu and P. Zhang (eds.), Number Theory and Related Area, Higher Education Press $\&$ International Press,
Beijing and Boston, 2013, pp. 149--197.

\bibitem{Su5} Z.-W. Sun, {\it On sums related to central binomial and trinomial coefficients},
in: M. B. Nathanson (ed.), Combinatorial and Additive Number Theory: CANT 2011 and 2012,
¡¡¡¡Springer Proc. in Math. $\&$ Stat., Vol. 101, Springer, New York, 2014, pp. 257--312.

\bibitem{vH} L. van Hamme, {\it Some conjectures concerning partial sums of generalized hypergeometric series},
in: $p$-adic Functional Analysis (Nijmegen, 1996), pp. 223--236, Lecture Notes in Pure and Appl. Math.,
Vol. 192, Dekker, 1997.


\end{thebibliography}
\end{document}